\documentclass{article}

  \usepackage{amsmath}
  \usepackage{paralist}
  \usepackage{graphics}
  \usepackage{epsfig}
  \usepackage[latin1]{inputenc}
  \usepackage[english]{babel}
  \usepackage{amssymb}
  \usepackage{hyperref}
  \usepackage{theorem}
  \usepackage{xypic}

\newtheorem{theorem}{Theorem}
\newtheorem{corollary}{Corollary}
\newtheorem{lemma}{Lemma}
\newtheorem{proposition}{Proposition}

{ \theorembodyfont{\rmfamily}
\newtheorem{definition}{Definition}
\newtheorem{remark}{Remark}
\newtheorem{example}{Example} }

\begin{document}

\title{A note on the relation between joint and differential invariants}
\author{David Bl\'azquez-Sanz \& Juan Sebasti\'an D\'iaz Arboleda}

\maketitle

\begin{abstract}
We discuss the general properties of the theory of joint invariants of a smooth
Lie group action in a manifold. Many
of the known results about differential invariants, including Lie's finiteness theorem,
have simpler versions in the context of joint invariants. We explore
the relation between joint and differential invariants, and  
we expose a general method that allow to compute differential invariants from
joint invariants.
\end{abstract}


\section{Introduction}

  Let us consider a connected Lie group $G$ and an smooth $G$-manifold $M$,
  $$\alpha\colon G\times M\to M, \quad (q,p) \mapsto \alpha(g,p) = gp.$$
For each $g\in G$ we denote by $\alpha_g$ the diffeomorphism that maps 
each point $p$ of $M$ to $gp$.   
The action of $G$ extends diagonally to 
all the Cartesian powers $M^k$ of $M$ by setting for each $k$-tuple 
$\bar p = (p_1,\ldots, p_k),$
$$g \bar p = (gp_1,\ldots,gp_k).$$
A smooth \emph{$k$-joint invariant} (or joint invariant of $k$ points) 
of the action of $G$ in $M$ defined in $W \subseteq M^k$
is a function $I\in\mathcal C^{\infty}_{M^k}(W)$ such that 
for each $k$-tuple  $\bar p\in W$ and each $g\in \mathcal G$ 
such that $g\bar p\in W$ we have $I(\bar p) = I(g \bar p)$.

Joint invariants appear frequently in classical geometry (see \cite{Olver1999, Olver2001}). 
Some well know
examples are: affine ratio in affine geometry; anharmonic 
ratio in projective geometry; distance of two points, area of a triangle, 
volume of a tetrahedron in Euclidean geometry. 
There is also a natural relation between joint invariants and the
geometry of differential equations. The joint invariants
of the Lie-Vessiot-Guldberg algebra of an ordinary differential equation
give rise to non-linear superposition formulas. 
Such relation has been studied
recently by several authors, \cite{Blazquez2010, Carinena2000, Carinena2011}. 
The most know example is the conservation 
of the anharmonic ratio by the Riccati equation
that gives rise to the non-linear superposition of three different 
solutions.

On the other hand, the action of the group $G$ in $M$, also
prolongs to the Weil near-point\footnote{Some authors prefer to define differential
invariants in the Jet bundles, that are algebraic quotients of the near-point bundles.
We prefer here the formalism of Weil near-point bundles because it adapt better to our
computations. A comparison between those notions, and detailed explanations about
the construction of Jet and Weil bundles can be found in \cite{Munoz2000}.}
bundles $T^{m,r} M$, of Taylor developments
of order $r$ of maps $(\mathbb R^m,0)\to M$. 
Let
$s\colon (\mathbb R^m, 0)\to M$ is a germ smooth map, and $j^r_0(s)\in T^{m,r}M$ his 
Taylor development of order $r$ at $0$. For each $g\in G$ with 
we can define, 
$$g(j_0^r(s)) = j_0^r(\alpha_g\circ s).$$ 
A smooth \emph{differential invariant} of order $r$ and rank $m$ of the action of $G$ in
$M$ is a function $I$ defined in an open subset $W\subseteq T^{m,r}M$ that is invariant
by the induced action of $G$ in $T^{m,r}M$. 
Differential invariants play an interesting role in the 
analysis of differential equations, the non-linear 
differential Galois theory, and $G$-structures.

There are some known natural relations between joint and
differential invariants.  For instance, let us consider 
a diffeomorphism $f$ of the projective line and $\{\sigma_t\}_{t\in \mathbb R}$ a
monoparametric group of transformations. 
It is well know (see for instance \cite{Ovsienko2005}, page 10) that the Schwartzian 
derivative can be seen as the infinitesimal deformation of the anharmonic ratio in the following way.
Let $x$ be any point, then
$$\frac{(f(x) - f(\sigma_{2\varepsilon}x))(f(\sigma_{\varepsilon})x - f(\sigma_{3\varepsilon}x))}
{(f(x) - f(\sigma_{\varepsilon}x))(f(\sigma_{2\varepsilon}x) - f(\sigma_{3\varepsilon}x))} =$$ 
$$\frac{(x - \sigma_{2\varepsilon}x)(\sigma_{\varepsilon}x - \sigma_{3\varepsilon}x)}
{(x - \sigma_{\varepsilon}x)(\sigma_{2\varepsilon}x - \sigma_{3\varepsilon}x)} 
+ \frac{3f''(x)^2 - 2f'(x)f'''(x)}{f'(x)^2}\varepsilon^2 + o(\varepsilon^3) $$

In this article we study the general theory of joint invariants, 
and explore the relation between joint and differential invariants. We show
that there is a general mechanism of derivation of joint invariants
that yields differential invariants.

\section{Joint invariants}

\subsection{Sheaf of local invariants}

For each $p\in M$ we denote by $Gp$ the orbit of $p$,
$$Gp = \{gp\in M\;|\, g\in G\}.$$
Let $U\subseteq M$ be an open subset. 
An smooth function $I\in\mathcal C_M^\infty(U)$ is called an
invariant if for each $p\in U$ and each $g\in G$ such that $gp\in U$ we have
$I(p) = I(g(p))$. In other words, for each $p\in U$ we have that $I|_{Gp\cap U}$
is a constant function. The functor that assigns to each open subset $U$ the
set of the invariants defined in $U$ is a presheaf. Its associated sheaf $\mathcal A^G$
is called the \emph{sheaf of local invariants} of the action of $G$ in $M$. We have thus,
$$\mathcal A^G(U) = \{I\in\mathcal C_M^\infty(U) \;|\; 
\forall p\in U\, I \mbox{ is locally constant in }Gp\cap U\}.$$
If $U$ is union of orbits of the action of the connected
Lie group $G$ then $\mathcal A^G(U)$ is exactly the ring of invariants of the action of $G$ in $U$.
For a general open subset $U$ the ring of local
invariants $\mathcal A^G(U)$ may be bigger that the ring of invariants.

\begin{example}\label{Proj1}
Let us consider the group ${\rm PGL}(2,\mathbb R)$ of classes modulo scalars of non-degenerated 
$2\times 2$ matrices. The M\"obius action of ${\rm PGL}(2,\mathbb R)$ in $\mathbb{RP}_1$
is given in the canonical affine coordinate $x$ in $\mathbb{RP}_1$ by the formula:
$$
\left[
\begin{array}{cc} 
a & b \\ c & d
\end{array}
\right]x = \frac{ax+b}{cx+d}.
$$
Let us consider the subgroups ${\rm Mov}(1,\mathbb R) \subset {\rm Aff}(1, \mathbb R) \subset 
{\rm PGL}(2,\mathbb R)$ where the respective inclusions are explicitly given as:
$$\left\{ \left[
\begin{array}{cc} 
1 & \lambda \\ 0 & 1
\end{array}
\right] \colon \lambda\in \mathbb R
\right\} \subset 
\left\{ \left[
\begin{array}{cc} 
\lambda & \mu \\ 0 & 1
\end{array}
\right] \colon \lambda,\mu \in \mathbb R, \lambda\neq 0 
\right\}
\subset 
{\rm PGL}(2,\mathbb R).
$$
We have:
\begin{itemize}
\item[(a)] The oriented distance $d^+(x_1,x_2) = x_2-x_1$ is a $2$-joint invariant in the open
subset $\mathbb R^2\subset (\mathbb{RP}_1)^2$ of the action of ${\rm Mov}(1, \mathbb R)$ in $\mathbb{RP}_1$.
\item[(b)] The affine ratio $[x_1,x_2;x_3] = \frac{x_3-x_1}{x_2-x_1}$ is a $3$-joint 
invariant in the open subset $\mathbb R^3\setminus \{x_2=x_1\} 
\subset (\mathbb{RP}_1)^3$ of the action of 
${\rm Aff}(1)$ in $\mathbb{RP}_1$.
\item[(c)] The anharmonic ratio 
$[x_1,x_2,x_3;x_4] = 
\frac{(x_1-x_3)(x_2-x_4)}{(x_1-x_2)(x_3-x_4)}$ is a $4$-joint 
invariant in the open subset $(\mathbb{RP}_1)^4\setminus \{x_i=x_j\}_{i\neq j} 
\subset (\mathbb{RP}_1)^4$ of the action of 
${\rm PGL}(2)$ in $\mathbb{RP}_1$.
\end{itemize}
\end{example}

\begin{example}\label{Euclid1}
Let us consider the group ${\rm Mov}(n,\mathbb R)$ of Euclidean motions in $\mathbb R^n$.
Let us write $\bar p = (p_1,\ldots,p_k)$ and $p_i = (x_{1,i},\ldots,x_{n,i})$.
\begin{itemize}
\item[(a)] The square distance $d^2(p_1,p_2) = \sum_{i=1}^n \left(x_{i,1} - x_{i,2}\right)^2$
is a smooth $2$-joint invariant defined in $(\mathbb R^n)^2$. 
\item[(b)] The Euclidean area $A(p_1,p_2,p_3)$ of a triangle of vertices $p_1$, $p_2$, $p_3$
is a smooth $3$-joint invariant defined in the open subset of $3$-tuples
of points in general position. Let us note that this set is dense if $n\geq 2$ and empty
for $n<2$.
\item[(c)] The Euclidean volume $V(p_1,p_2,p_3,p_4)$ of tetrahedron of vertices $p_1$, $p_2$, $p_3$,
$p_4$ is a smooth $4$-joint invariant defined in the open subset of 
$4$-tuples of points in general position. Let us note that this set is dense if $n\geq 3$ and empty
for $n<3$.
\item[(d)] The $k$-dimensional Euclidean volume $V_k(p_1,\ldots,p_{k+1})$ of 
a $k$-simplex of vertices $p_1$,$\ldots$, $p_{k+1}$ is a smooth
$(k+1)$-joint invariant defined in the open subset of 
$(k+1)$-tuples of points in general position. Let us note this set is dense if 
$n> k$ and empty
for $n \leq k $.
\end{itemize}
\end{example}

\subsubsection*{Infinitesimal generators}

Let us denote by $\mathcal G$ the Lie algebra of of right invariant vector fields in $G$.
There is a natural Lie algebra morphism from $\mathcal G$
into the into the Lie algebra $\mathfrak X(M)$ of
vector fields in $M$:
$${\bf ig}\colon \mathcal G \to \mathfrak X(M),  
\quad X \mapsto {\bf ig}(X),\quad  {\bf ig}(X)_p = \left.\frac{d(e^{tX}p)}{dt}\right|_{t=0}. 
$$
We call infinitesimal generators of the action of $G$ in $M$ to the vector fields
in the image of ${\bf ig}$. The sheaf of local invariants $\mathcal A^G$ is easily characterized as
the sheaf of first integral of the Lie algebra of the infinitesimal generators in $M$. 
$$\mathcal A^G(U) = \{f\in\mathcal C^\infty(U), \;|\; \forall X\in\mathcal G\;\; {\bf ig}(X)f = 0\}.$$

The following result tell us that we can restrict our consideration to dense open
subsets of $M$.

\begin{lemma}\label{Lemma1}
Let $U\subseteq W \subseteq M$ be open subsets in $M$ such that $U$ is dense in $W$. 
Then if $I$ is a function in $W$ such that $I|_U$ is a local
invariant, then $I$ is  local invariant in $W$. 
\end{lemma}

\proofname. Let $A_1,\ldots,A_r$ be a basis of $\mathcal G$. Then, $I$
is a local invariant in $U$ if and only if ${\bf ig}(A_i)I = 0$ in $U$
for $i=1,\ldots,r$. By continuity we have that ${\bf ig}(A_i)I = 0$ in $W$
and thus $I$ is a local invariant in $W$. 
\hfill$\square$

\medskip

If is also useful to consider the distribution of vector fields spanned by the infinitesimal 
generators. We denote by $\mathcal L^G$ the distribution of vector fields that assigns to each 
point $p$
of $M$ the vector space $\mathcal L^G_p$,
 $$\mathcal L^G_p = {\rm span}\{{\bf ig}(X)_p \,|\, X\in \mathcal G\} \subseteq T_pM,$$
it is clear that $T_p(Gp) = \mathcal L^G_p$. The distribution $\mathcal L^G$ is,
by definition, stable by the Lie bracket. It determines
a unique Pfaff system $\Sigma^G$ such that $(\Sigma^G)^\perp = \mathcal L^G$. We review
the basic definitions about Pfaff systems and their first integrals in the next section.

\subsection{Some considerations about first integrals}

We denote by $\Omega^\bullet_M$ the sheaf of differential
exterior forms in $M$. Thus, 
$$\Omega^\bullet_M(U) = \bigoplus_{k=0}^{\dim (M)} \Omega^k_M(U).$$
For open $p\subset M$ the stalk $\Omega^\bullet_{M,p}$  
is endowed with the wedge product and the exterior differential
it is a free $\mathcal C^\infty_{M,p}$
module of rank $2^{\dim(M)}$.
Thus, $\Omega^\bullet_M$ is a sheaf of differential rings, and also
a locally free of finite rank sheaf of $\mathcal C_M^\infty$-modules.
Given a a set $S\subset \Omega^\bullet_{M,p}$ we will denote
by $(S)$ the ideal of $\Omega_{M,p}^\bullet$ spanned by $S$ and 
by $\{S\} = (S,dS)$ the differential ideal of $\Omega_{M,p}^\bullet$ spanned by $S$.
The same notation extends easily for a subsheaf $\mathcal S\subset \Omega^\bullet_M$.

\begin{definition}
An exterior differential system (e.d.s.) $\Sigma$ in $M$ is a sheaf $\Sigma\subset\Omega_M^\bullet$ 
of differential ideals without zero-order equations. That is, for each $p\in M$ it
satisfies $\Sigma_p\cap \Omega^0_{M,p} = \{0\}$.
\end{definition}

\begin{definition}
A differential exterior system $\Sigma$ in $M$ is a Pfaff system if
it is spanned, as a sheaf of differential ideals, by its sections of
degree one, that is for each $p\in M$, $\Sigma_p = \{\Sigma_p\cap \Omega^1_{M,p}\}$.   
\end{definition}

\begin{definition} Let $U$ be an open subset in $M$.  
A function $F\in\mathcal C^\infty(U)$ is called a first integral of 
$\Sigma$ if $dF\in \Sigma(U)$. The first integrals of $\Sigma$
form a sheaf of smooth functions that we denote by ${\bf int}(\Sigma)$.
\end{definition}

First integrals commute with smooth maps. If $\varphi\colon N\to M$
is smooth and $\Sigma$ is an e.d.s. in $M$, then ${\bf int}(\varphi^*\Sigma) = 
\varphi^*({\bf int}(\Sigma))$.  
Given a sheaf $\mathcal S\subset \mathcal C^\infty_M$ of rings smooth functions 
in $M$ we can 
differentiate it to obtain a Pfaff system $d\mathcal S$,
$$(d\mathcal S)_p = \{d(\mathcal S_p)\}\subset \Omega^\bullet_{M,p}.$$
It is clear that all sections of $\mathcal S$ are first integrals
of $d\mathcal S$,
$$\mathcal S \subseteq {\bf int}(d\mathcal S).$$
And, for a Pfaff system $\Sigma$, the differential of first integrals of
$\Sigma$ are by definition in $\Sigma$,
$$d({\bf int}(\Sigma)) \subseteq \Sigma.$$ 

\begin{definition}
A Pfaff system $\Sigma$ is called \emph{integrable} if it completely determined
by its sheaf of first integrals, $d({\bf int}(\Sigma)) = \Sigma.$   
\end{definition}

A Pfaff system $\Sigma$ has a associated several distribution of
vector fields $\Sigma^\perp$ and $\Sigma'$ called the orthogonal
and de characteristic distribution. For each $p\in M$ we have:
$$\Sigma^\perp_p = \{X_p\in T_pM \;|\; \forall \omega\in 
\Sigma_p\cap\Omega^1_p\;\;\omega_p(X_p) = 0  \}$$
$$\Sigma'_p = \{X_p\in T_pM \;|\; \forall \omega\in 
\Sigma_p\;{\bf i}_{X_p}\omega_p = 0  \}$$
It is clear that $\Sigma'\subseteq \Sigma^\perp$. The distribution
$\Sigma'$ stable by Lie brackets. 
 
\begin{definition}
A Pfaff system $\Sigma$ is called regular of rank $r$ if for each $p\in M$
there are $\omega_1,\ldots,\omega_r\in \Sigma_p$ such that:
\begin{itemize}
\item[(a)] $\Sigma_p = \{\omega_1,\ldots,\omega_r\}$.
\item[(b)] $\omega_1(p),\ldots,\omega_r(p)\in T^*_pM$ are linearly
independent. 
\end{itemize}
\end{definition}

It is clear that $\Sigma$ is a regular Pfaff system of rank $r$ if and only if
the distribution of vector fields $\Sigma^{\perp}$ is a regular sub-bundle 
of $TM$ of rank $\dim(M) - r$. Thus,
regular distributions of vector fields and regular Pfaff systems in $M$ are
in one-to-one correspondence: a $1$-form $\omega$ is a section of $\Sigma$
if and only if it vanish along $\Sigma^\perp$. 
 It is also interesting to remark that
for any regular submersion $\varphi\colon N \to M$ we have 
$\varphi^*(\Sigma)$ is a regular and: $$\varphi^*(\Sigma)^\perp = d\varphi^{-1}(\Sigma^{\perp}).$$ 
 The following result is well known:

\begin{theorem}[Frobenius]
Let $\Sigma$ be a regular Pfaff system in $M$ of rank $r$. The following are equivalent:
\begin{enumerate}
\item[(a)] $\Sigma$ is integrable.
\item[(b)] $\Sigma^\perp$ is stable by Lie bracket.
\item[(c)] $\Sigma^\perp = \Sigma'$.
\item[(d)] $\Sigma$ is spanned as a sheaf of ideals its homogeneous component 
of degree one:
$\Sigma = (\Sigma\cap\Omega^1_M).$
\item[(e)] Through each point $p\in M$ it passes a submanifold $S$ of dimension 
$\dim M - r$ such that for each $q\in S$, $T_qS = \Sigma^\perp_q$.
\item[(f)] For each point $p\in M$ there are $r$ functionally
independent first integrals of $\Sigma$ in $\mathcal C^\infty_{M,p}$.
\end{enumerate}
\end{theorem}

A system of $r$ functionally independent 
first integrals of $\Sigma$, defined on an open subset $U$, is called a 
\emph{complete system} of first integrals of $\Sigma$ in $U$. Any other first 
integral of $\Sigma$ is locally functionally dependent of them.

\subsection{Functional dependence}

Let us consider $\mathcal R$ and $\mathcal S$ two sheaves of smooth functions
on $M$. We can consider the sheaf of functions 
that locally depend functionally of sections of $\mathcal R$ and $\mathcal S$.
The theorem of functional dependence says that 
a germ of smooth function $f\in\mathcal C^\infty_{M,p}$ is function 
of germs $f_1,\ldots,f_k\in \mathcal R_p$, $h_1,\ldots,h_s\in \mathcal S_p$
if there are $\lambda_1,\ldots,\lambda_{k+s} \in \mathcal C^\infty_{M,p}$ such that
$$df = \sum_{i=1}^k \lambda_idf_i + \sum_{j=1}^s \lambda_{k+j}dh_j.$$
Thus, a function $f\in \mathcal C^\infty(U)$ 
is locally functionally dependent of those in $\mathcal R(U)$ and $\mathcal S(U)$
if and only if $f$ is a simultaneously first integral of the Pfaff systems
$d\mathcal R$ and $d\mathcal S$. 
This consideration allow us to give a simpler definition:

\begin{definition}
The sheaf $\mathcal R\odot \mathcal S$ of functions that locally depend functionally
of those in $\mathcal R$ and $\mathcal S$ is:
$$\mathcal R\odot \mathcal S 
= {\bf int}(\{d\mathcal S,d \mathcal R\})$$
\end{definition}

The following result list some direct consequences of the definitions and does
not need a proof.

\begin{lemma}\label{Lemma2}
Let $\Sigma_1$ and $\Sigma_2$ be two regular Pfaff systems in $M$ such that
$\{\Sigma_1,\Sigma_2\}$ is also regular, 
and let $\mathcal R_1$, $\mathcal R_2$ be their
respective sheaves of first integrals. 
\begin{itemize} 
\item[(a)] $\{\Sigma_1,\Sigma_2\}$ is integrable and 
$\mathcal R_1\odot\mathcal R_2$ is its sheaf of first integrals.
\item[(b)] $f\in \mathcal C_M^\infty(U)$ for $U$ open subset of $M$ is a first integral
of $\{\Sigma_1,\Sigma_2\}$ if and only if $df$ vanish on the 
distribution of vector fields
$\Sigma_1^\perp \cap \Sigma_2^\perp$.
\end{itemize}
\end{lemma}

\subsection{Lifting of joint invariants}

From now on let us fix some notation about the sheaves of local joint invariants.
Let us denote $\mathcal A^{G,k} \subset \mathcal C^\infty_{M^k}$ the sheaf
of local $k$-joint invariants, $\mathcal L^{G,k}$ the distribution of vector fields
spanned by the infinitesimal generators of the action of $G$ in $M^k$ and $\Sigma^{G,k}$ the
Pfaff system generated by $1$-forms in open sets of $M^k$ vanishing on $\mathcal L^{G,k}$.

By a \emph{generating system of (local) $k$-joint invariants} we mean a 
system of first integrals 
of $\Sigma^{G,k}$. If we work in an open subset of $M^k$ in which $\mathcal L^{G,k}$
is a regular distribution of vector fields of the same dimension than $G$, the
Frobenious theorems ensures that all 
generating systems of (local) $k$-joint invariants consist of $k \dim(M)-\dim(G)$
functionally independent invariants.  

\begin{lemma}\label{Lemma4}
Assume that $M^k$ contains an open subset $W_k$ in which $\Sigma^{G,k}$ is regular.
Let $U\subseteq M^k$ be an open subset. 
Any smooth function $I\in\mathcal C^\infty(U)$ with $U\subset M^k$ is a local
$k$-joint invariant if and only if $I|_{W_k\cap U}$ is a first integral of $\Sigma^{G,k}$. 
\end{lemma}

\proofname. If $\Sigma^{G,k}$ is regular in $W_k$, a first integral in $W_k$ is a 
funtion wich is locally constant along the orbits, and thus, a local $k$-joint invariant.
Then, by Lemma \ref{Lemma1}, we conclude.
\hfill$\square$

\medskip

We also denote by $\pi_i^k$ the projection from $M^k$ to $M^{k-1}$ that consists in
dropping the $i$-th component:
$$\pi^k_i\colon M^k\to M^{k-1}, \quad (p_1,\ldots,p_k)\mapsto (p_1\ldots,p_{i-1},p_{i+1},\ldots,p_k).$$
These projections $\pi_i^k$ are $G$-equivariant, and thus: 
$$\pi_i^{k*}(\mathcal A^{G,k-1}) \subseteq \mathcal A^{G,k}.$$
This inclusion simply means that each $(k-1)$-joint invariant of can
be seen as a $k$-joint invariant in $k$ different ways by dropping one of the arguments.

Let $W\subset M^k$ be an open subset. Let us consider a 
function $I\in \mathcal C^{\infty}_{M^k}(W)$ that depends functionally 
on local joint invariants of $k-1$ points:
$$I \in (\pi_1^{k*}\mathcal A^{G,k-1}\odot \ldots \odot \pi_k^{k*}\mathcal A^{G,k-1})(W),$$
then it is clear that it is a local joint invariant of $k$ points, thus:
$$\pi_1^{k*}\mathcal A^{G,k-1} \odot \ldots \odot \pi_k^{k*}\mathcal A^{G,k-1} 
\subseteq \mathcal A^{G,k}.$$
We can ask if any local joint invariant of $k$ points can be expressed as a function of
local joint invariants of less that $k$ points. 

\begin{example}
Let us consider the action of ${\rm PGL}(2,\mathbb R)$ ant its subgroups in $\mathbb{RP}_1$ 
discussed in Example \ref{Proj1}. Since $\mathbb{RP}_1$ is $1$-dimensional, and by two
different liftings $\pi_i^k\colon\mathbb{RP}_1^k\to \mathbb{RP}_1^{k-1}$ we will always
give functionally independent joint invariants, we can conclude:
\begin{itemize}
\item[(a)] Any $k$-joint invariant of 
${\rm Mov}(1,\mathbb R)$ with $k\geq 2$ is functionally
dependent of the liftings of the oriented distance.
\item[(b)] Any $k$-joint invariant of 
${\rm Aff}(1,\mathbb R)$ with $k\geq 3$ is functionally
dependent of the liftings of the affine ratio
\item[(a)] Any $k$-joint invariant of ${\rm PGL}(2,\mathbb R)$ with $k \geq 4$ 
is functionally dependent of the liftings of the anharmonic ratio.
\end{itemize}   
\end{example}

\begin{example}
Let us discuss the functional dependence of the joint invariants shown 
in Example \ref{Euclid1}. Let $\bar p = (p_1,\ldots,p_k)$
be a tuple in $(\mathbb R^n)^k$. By iterating $k-2$ times the liftings from 
$(\mathbb R^n)^2$ up to $(\mathbb R^n)^k$ we obtain ${k \choose 2}$ functions
called \emph{mutual distances}:
$$d_{ij}(\bar p) = d(p_i,p_j).$$
For $k\leq n$ the mutual distances are functionally independent in a dense open subset
of $(\mathbb R^n)^k$. For $k=n$ we have that the rank 
of $(\Sigma^{G,n})^\perp$
coincides with the dimension of ${\rm Mov}(\mathbb R,n)$ which is ${n+1 \choose 2 }$.
In this case the mutual distances form a complete system of ${n\choose 2}$ first integrals,
the dimension of the manifold of configurations of $n$ points is  
$(\mathbb R^n)^n$ is $n^2 = {n\choose 2} + {n+1 \choose 2}$, the sum of the 
number of independent $n$-joint invariants and the dimension of the group. 
The functional dependence relations between the mutual distances and the 
the rest of the $k$-joint invariants discussed in Example \ref{Euclid1} 
are well known (see \cite{Gritzmann1994} sec. 3.6.1) theorems of
Euclidean geometry, namely:
\begin{itemize}
\item[(a)] The area of a triangle $A(p_1,p_2,p_3)$ is a function of $d_{12}$, $d_{23}$,
and $d_{23}$ by means of the Heron formula:
$$A =$$ 
$$\frac{1}{4}\sqrt{(d_{12}+ d_{13} + d_{23})(-d_{12}+d_{13}+d_{23})
(d_{12}-d_{13}+d_{23})(d_{12}+d_{13}-d_{23})}.$$
We shall also remark that the classical theorems of trigonometry, the sinus and cosinus
theorems, are also relations of functional dependence between $3$-joint invariants
of the Euclidean plane geometry.
\item[(b)] The volumen of the tetrahedron $V(p_1,p_2,p_3,p_4)$ is a function
of the functions $d_{ij}$ for $i,j=1,\ldots,4$ by means of thr Cayley-Menger determinant:
$$V = \frac{1}{12\sqrt{2}}\sqrt{
\left|
\begin{array}{ccccc}
0 & 1 & 1& 1 & 1 \\
1 & 0 & d_{12}^2 & d_{13}^2 & d_{14}^2 \\
1 & d_{12}^2 & 0 & d_{23}^2 & d_{24}^2 \\
1 & d_{13}^2 & d_{23}^2 & 0 & d_{34}^2 \\
1 & d_{14}^2 & d_{24}^2 & d_{34}^2 & 0 
\end{array}
\right|
}$$
\item[(c)] The $k$-dimensional volume satisfies the equation:
$$2^k (k!)^2 V_k^2 = (-1)^{k+1}\det(A),$$
where $C$ is the $(k+1)\times (k+1)$ matrix whose elements $c_{ij}$ are:
$$
c_{ij} = \begin{cases} 
0  \quad\mbox{if}\quad i=j,\\
1  \quad\mbox{if}\quad i\neq j=1 \mbox{ or } j\neq i = 1, \\
d_{i-1,j-1}^2 \quad \mbox{in any other case.} \\
\end{cases} 
$$
\end{itemize}  
\end{example}

\subsection{Finiteness theorem}

For each $p$ in $M$ we denote by ${\rm Est}(p)$ the stabilizer group of $p$, 
$${\rm Est}(p) = \{g\in G \:|\: gp = p\}$$
and by ${\rm est}(p)$ its Lie algebra. 
With respect to the action of $G$ on $k$-tuples, it is useful to note that
for $\bar p = (p_1,\ldots,p_k)$,
$${\rm Est}(\bar p) = \bigcap_{i=1}^k {\rm Est}(p_i),\quad
{\rm est}(\bar p) = \bigcap_{i=1}^k {\rm est}(p_i).$$

The action of $G$ in $M$ is called \emph{locally faithful} 
if for each open subet $U\in M$ the map:
$${\bf ig}|_U\colon \mathcal G \to \mathfrak X(U), \quad A\mapsto {\bf ig}(A)|_U$$ 
is inyective. It is equivalent to say that for each $A\in \mathcal G$ the 
set of zeroes of ${\bf ig}(A)$ has
empty interior. It is clear that an algebraic or analytic faithful action
is also locally faithful, but an smooth faithful action can be not locally
faithful.

A point $p\in M$ is said \emph{locally regular} if ${\rm est}(p) = \{0\}$. The set of
locally regular points is an open subset of $M$. We will say that the action of $G$ in $M$
is \emph{generically locally regular} 
if the set of locally regular points is dense in $M$. 

It is clear that if the action of $G$ in $M^k$ is generically locally regular, 
then the action of $G$ in $M^{k+1}$ is also generically locally regular. If there
is an smaller natural number $k_0$ such that the action of $G$ in $M^{k_0}$ is
generically locally regular, we say that $k_0$ is the \emph{local rank} of
the action of $G$ in $M$. For $k>k_0$ the Pfaff system $\Sigma^{k,G}$, restricted
to the dense open subset of locally regular points, is a regular integrable system.

\begin{lemma}\label{Lemma3}
If the action of $G$ in $M$ is locally faithful then the action
of $G$ in $M$ has a finite local rank $k_0$. It verifies:
$$ \frac{\dim G}{\dim M} \geq k_0 \geq \dim G.$$ 
\end{lemma}

\proofname. 
Let $\bar p = (p_1,\ldots,p_k)\in M^k$. The vector space ${\rm est}(\bar p)$
is by definition the kernel of a linear map 
from $\mathcal G$ to $T_{\bar p}(M^k)$. The 
$k$-tuple $\bar p$ is regular if and only if this morphism is injective. This
is only possible if the dimension of the target space is bigger than the dimension 
of the source space. This gives us the left side of the inequality of the statement.

Let us see that the set of locally regular $k$-tuples, with $k\geq\dim G$ is 
dense. Let us consider $\bar p = (p_1,\ldots,p_k)\in M^k$, and $U$ any
neighbourhood of $\bar p$, without loss of generality we can assume that
$U = U_1,\ldots, U_k$ where each $U_i$ is a neighbourhood of $p_i$. Since
the action is locally faithful there is $p_1^*$ in $U_1$ such that the inclusion
${\rm est}(p_1^*) \subset \mathcal G$ is strict. Let us denote this space as $E_1$.
Then for each $i = 1,\ldots, k-1$ we set the space $E_k$ in the following manner:
\begin{itemize}
\item[(a)] If $E_{i} = \{0\}$ we set  $p_{i+1}^* = p_{i+1}$ and $E_{i+1} = \{0\}$.
\item[(b)] If If $E_{i}\neq \{0\}$ we take a non-zero vector $A_i\in E_i$. 
By hypothesis the set of zeroes of ${\bf ig}(A_i)$
is of empty interior in $U_{i+1}$ then, we set $p_{i+1}^*$ a point in $U_2$ such that 
${\bf ig}(A_i)(p_{i+1}^*) \neq 0$ and $E_{i+1} = {\rm est}(p_1^*,\ldots, p_{i+1}^*)$.
In this case we have $A_i\not\in E_{i+1}$, and thus an strict inclusion
$E_i\supset E_{i+1}$.
\end{itemize}
By the above argument, we have defined a descending chain of vector spaces:
$$\mathcal G \supseteq E_1 \supseteq \ldots \supseteq E_k,$$
where each inclusion is strict until the chain stabilizes at $0$. It follows
that for $k\geq \dim G$ then $E_k = 0$, thus:
$${\rm est}(p_1^*,\ldots,p_k^*) = E_k = 0,$$
and $\bar p^* = (p_1^*,\ldots,p_k^*)$ is a regular $k$-tuple un $U$. 
We have proven that for $k\geq \dim G$ the action of $G$ in $M^k$ is
generically regular. \hfill$\square$

\medskip

\begin{remark}
The left inequality 
of the statement is treated by S. Lie in the study of
superposition laws, in such context it is know as Lie inequality. 
This local rank coincides is the number of known solutions
necessary to describe the general solution in a o.d.e. admitting non linear
superposition laws. The reader may consult 
\cite{Blazquez2010, Carinena2000, Carinena2011} for further explanations.
\end{remark}

\begin{proposition}\label{prop_fin}
Assume that $\Sigma^{G,k-1}$ is a regular integrable Pfaff system in some dense 
open subset $W_{k-1}$ of $M^{k+1}$ and that there is an dense open subset 
$W_k$ in $M^k$ such that
for each $k$-tuple $\bar p\in W_k$ we have:
\begin{equation}\label{cond_lineal}
\bigcap_{i=1}^k \left(\mathcal L^{G,k}_{\bar p} \oplus \ker(d_{\bar p}\pi_i^k)\right) = \mathcal L^{G,k}_{\bar p},
\end{equation}
then each local joint invariant of $k$ points can be locally expressed as a function of local
joint invariants of $k-1$,
$$\mathcal A^{G,k} = 
\pi_1^{k*}\mathcal A^{G,k-1} \odot \ldots \odot \pi_k^{k*}\mathcal A^{G,k-1}. 
$$ 
\end{proposition}

\proofname. 
We have that 
$$\left(\pi_i^{k*}(\Sigma^{G,k-1})\right)^\perp = 
\mathcal L^{G,k} \oplus \ker(d\pi_i^k),$$
thus, the condition \eqref{cond_lineal} is equivalent to: 
$$\{\pi_1^{k*}(\Sigma^{G,k-1}),\pi_2^{k*}(\Sigma^{G,k-1}), 
\ldots, \pi_i^{k*}(\Sigma^{G,k-1})\} = \Sigma^{G,k}$$
and then by Lemmas \ref{Lemma1} and \ref{Lemma2} we finish.  
\hfill$\square$

\medskip

After integration, condition \eqref{cond_lineal} have the following geometrical meaning. 
Let $U\subset W_k$ an small enough subset of an orbit in $M^k$. Then, $U$ can be
recovered as the intersection of its projections onto $M^{k-1}$:
$$U = \bigcap_{i=1}^k (\pi^{k}_i)^{-1}(\pi_i^k(U)).$$

We will see that condition \eqref{cond_lineal} is always satisfied for locally faithful actions.
If fact, for $k$ big enough we only need to take into account two factors of the intersection,
as the following elementary lemma of linear algebra shows.

\begin{lemma}\label{linear_inter}
Let $A,B,C$ be real vector subspaces of a real vector space $E$.  
If $B\cap C = \{0\}$ and 
$A\cap (B\oplus C) = \{0\}$ then $(A\oplus B)\cap (A \oplus C) = A$.
\end{lemma}

\proofname. We prove only the non trivial inclusion. 
Let $v$ be in $(A\oplus B)\cap (A \oplus C)$, then there are
decompositions $v= a_1+b = a_2 + c$ with $a_1,a_2\in A$, $b\in B$ and $c\in C$.
Let us take $w = a_1-a_2 = c - b$, it is clear $w \in A\cap(B\oplus C)$ and then
$w = 0$, it yields $a_1 = a_2$ and hence $b = c = 0$ so that $v = a_1$ what
proves $v\in A$.\hfill$\square$

\begin{theorem}[Finiteness theorem]
Assume that the action of $G$ in $M$ is locally faithful, and let $k_0$ be its
local rank. There is a $k_1$ such that any local $k$-joint invariant in $M$
points with $k>k_1$ is functionally dependent of the liftings of local $k_1$-joint invariants. 
This number $k_1$ satisfies,
$$k_1\leq k_0 + 2.$$ 
\end{theorem}

\proofname. Let us prove that for each $k > k_0+2$ the local $k$-joint invariants
are functionally dependent of local $(k-1)$-joint invariant. By Lemmas
\ref{Lemma4} and \ref{Lemma3} we consider
a dense open subset $W_k$ in $M^k$ such that the action of $G$ in $M^k$
is locally regular each on $\pi_1^{k-1}(\pi_1^k(W_k))\subset M^{k-2}$. For each $\bar p$
in $W_k$ we have that the spaces $\mathcal L_{\bar p}^{G,k}$, 
$\ker(d_{\bar p} \pi_1^k)$ and $\ker(d_{\bar p} \pi_2^k)$ satisfy the hypothesis
of Lemma \ref{linear_inter}. Thus, we can apply Proposition \ref{prop_fin}, 
finishing the proof. \hfill$\square$

\section{Weil near-points bundles}

  In order to compute differential invariants, we will use the formalism of Weil bundles.
This is an approach to differential geometry proposed by Andre Weil in 1953 \cite{Weil1953}, 
who introduced the notion of infinitesimally near points. 
This formalism is an alternative to the better
known notion of Jet bundles of Ehresmann. It has
been developed independently by Shurygin \cite{Shurygin1989},  
Kol\'a\v{r}, Michor, Slov\'ak \cite{Kolar1993}, and Mu\~noz, 
Rodr\'iguez, and Muriel \cite{Munoz2000}. The reader interested in the proofs
of the statements and the details of the theory presented in this section is encouraged
to consult the former references.


\subsection{Infinitesimally near points}

 The spectral representation theorem says that the there is a canonical 
bijection from $M$ to ${\rm Hom}_{\mathbb R-alg}(\mathcal C^\infty(M),\mathbb R)$. To each
point $p \in M$ it corresponds the valuation morphism:
$$p\colon \mathcal C^\infty(M) \to \mathbb R, \quad f\mapsto p(f) := f(p).$$ 

\begin{definition}
A Weil algebra $A$ is a finite dimensional local $\mathbb R$-algebra with maximal
ideal $\mathfrak m_A$ whose quotient field $A/\mathfrak m_A$ is the field $\mathbb R$ of real numbers.
The quotient morphism $\omega_A\colon A\to\mathbb R$ is called the valuation morphism of $A$.
\end{definition}

In the applications we mainly use a specific kind of Weil algebras, the algebras of truncated
Taylor series of order $r$ in $m$ variables,
$$\mathbb R[[\varepsilon]]_{m,r} = \mathbb R[[\varepsilon]]/(\varepsilon)^{r+1},\quad\quad
\varepsilon = (\varepsilon_1,\ldots,\varepsilon_m).$$


Let us consider a Weil algebra $A$.
We define the bundle of \emph{$A$-near-points}\footnote{$\mathbb R[[\varepsilon]]_{m,r}$-near-points
are simply termed $(m,r)$-near-points.} in $M$
as $M(A) = {\rm Hom}_{\mathbb R-alg}(\mathcal C^\infty(M),A)$. By composition with
the valuation $\omega_A$ we have a canonical projection $\pi_A\colon M(A)\to M$,
which is an smooth bundle. For each open coordinate subset $U\subset M$
with coordinate functions,
$$ \bar x \colon U \mapsto W\subset \mathbb R^n$$
we have that $\pi^{-1}_A(U) = U(A)$ is endowed with coordinates with values in $A$,
$$\bar x^A\colon U(A)\to A^n,$$
and that system of coordinates identify $U(A)$ with the open subset of $A^n$
consisting in $n$-tuples $\bar a = (a_1,\ldots, a_n)$ such that $(\omega(a_1),\ldots,\omega(a_n))$
is in $W$.

We also consider the open sub-bundle $M(A)_{\rm prop}$ consisting
in the $A$-near points that are surjective as algebra morphisms. We have
some interesting and self-explanatory examples:

\begin{enumerate}
\item[(a)] The bundle $M(\mathbb R[[\varepsilon]]_{1,1})\to M$ is the tangent bundle $TM\to M$.
A vector $\vec X_p\in T_pX$ is seen as an $\mathbb R$-algebra morphism in the following way:
$$f \mapsto f(p) + \varepsilon \vec X_pf.$$
\item[(b)] The bundle $M(\mathbb R[[\varepsilon]]_{m,1})\to M$ is the bundle of frames $F(M)\to M$.
A frame $(\vec X_{1p},\ldots, \vec X_{np})$ is seen as an $\mathbb R$-algebra morphism if
the following way:
$$f \mapsto f(p) + \varepsilon_1 \vec X_{1p}f + \ldots \varepsilon_n \vec X_{np}f.$$
\item[(c)] The bundle $M(\mathbb R[[\varepsilon]]_{1,r})\to M$ is the tangent bundle of order $r$,
that we denote by $T^r M\to M$. 
An element $j^r_0\gamma$ of $T^rM$ is the Taylor development of order $r$ at $0$
of a curve $\gamma\colon(\mathbb R^m,0)\to M$. If the expression in coordinates of $\gamma$ is:
$$\gamma(\varepsilon) = (x_1(\varepsilon),\ldots,x_n(\varepsilon)),$$
where,
$$\begin{cases} 
x_1(\varepsilon) = \gamma_{1} + \gamma_{1}'\varepsilon + \ldots + 
\frac{\gamma_1^{(r)}}{r!}\varepsilon^r + 
o(\varepsilon^{r+1}),  \\ \vdots \\ 
x_n(\varepsilon) =  \gamma_{n} + \gamma_{n}'\varepsilon + \ldots + 
\frac{\gamma_n^{(r)}}{r!}\varepsilon^r +
o(\varepsilon^{r+1}), 
\end{cases}$$ 
then, $j^r_0\gamma$ is the morphism:
$$j^r_0\gamma\colon x_i \mapsto 
\gamma_{i} + \gamma_{i}'\varepsilon + \ldots + \frac{\gamma_i^{(r)}}{r!}\varepsilon^r
\in \mathbb R[[\varepsilon]]_{1,r}.$$
By convention we say that the coefficient $\gamma_i^{(j)}$ of the near-point
is its $x_i^{(j)}$ coordinate. Thus, if $x_1,\ldots,x_n$ is a system of coordinates
in $U\subseteq M$ then $x_1,\ldots,x_n$,$x_1',\ldots,x_n'$,$\ldots$,
$x_1^{(r)},\ldots,x_n^{(r)}$ is a system of coordinates in $T^rU\subseteq T^rM$. 
\item[(d)] The bundle $M(\mathbb R[[\varepsilon]]_{m,r})\to M$ is the generalized 
tangent bundle of order rank $m$ and order $r$,
that we denote by $T^{m,r} M\to M$. 
An element $j^r_0\varphi$ of $T^{m,r}M$ is the Taylor development of order $r$ at $0$
of a smooth map $\varphi\colon(\mathbb R^m,0)\to M$. If the expression in coordinates of $\gamma$ is:
$$\varphi(\varepsilon_1,\ldots,\varepsilon_m) = 
(x_1(\varepsilon_1,\ldots,\varepsilon_m),\ldots,x_n(\varepsilon_1,\ldots,\varepsilon_m)),$$
where,
$$\begin{cases} 
x_1(\varepsilon_1,\ldots,\varepsilon_m) = \varphi_{1,0} + \sum_{1\leq|\alpha|\leq r } 
\varphi_{1,\alpha}\frac{\varepsilon_1^{\alpha_1}\cdots \varepsilon_m^{\alpha_m}}
{\alpha_1!\cdots\alpha_m!}
+ o(\varepsilon^{r+1}),  \\ \vdots \\ 
x_n(\varepsilon_1,\ldots,\varepsilon_m) =  \varphi_{n,0} + \sum_{1\leq|\alpha|\leq r } 
\varphi_{n,\alpha}\frac{\varepsilon_1^{\alpha_1}\cdots \varepsilon_m^{\alpha_m}}
{\alpha_1!\cdots\alpha_m!} 
+ o(\varepsilon^{r+1}), 
\end{cases}$$ 
then, $j^r_0\varphi$ is the morphism:
$$j^r_0\varphi\colon x_i \mapsto 
\varphi_{i,0} + \sum_{1\leq|\alpha|\leq r } 
\varphi_{i,\alpha}\frac{\varepsilon_1^{\alpha_1}\cdots \varepsilon_m^{\alpha_m}}
{\alpha_1!\cdots\alpha_m!} 
\in \mathbb R[[\varepsilon]]_{1,r}.$$
By convention we say that the coefficient $\varphi_{i,alpha}$ of the near-point
is its $x_{i,\alpha}$ coordinate. Thus, if $x_1,\ldots,x_n$ is a system of coordinates
in $U\subseteq M$ then $x_1,\ldots,x_n$ and $x_{i,\alpha}$ with $i = 1,\ldots, n$,
and $\alpha = (\alpha_1,\ldots,\alpha_m)$ with $1\leq |\alpha| \leq r$ form a
a system of coordinates in $T^{m,r}U\subseteq T^{m,r}M$. 
\end{enumerate}

\subsection{Multi-type near points}


\begin{definition}
A \emph{multi-Weil algebra} $A$ is a finite direct product of Weil algebras. 
We will say that a point of 
$$M(A) = {\rm Hom}_{\mathbb R-alg}(\mathcal C^\infty(M), A)$$
is a multi-near-point of type $A$. 
\end{definition}

Multi-Weil algebras are not usually considered in the theory of near-points. 
However they
are useful in order to relate joint and differential invariants. 
Given a multi-Weil algebra, 
$$A = A_1\times \ldots \times A_k$$
the number $k$ of factors is called the \emph{multiplicity} of $A$. Note
that ${\rm Spec}(A)$ consist in $k$-maximal ideals $\{\mathfrak m_1, \ldots, \mathfrak m_k\}$ 
and there are $k$ valuation morphisms $\omega_{Ai}\colon A\to \mathbb R$. We donote by $\omega_A$
the map $(\omega_{A1}\ldots,\omega_{Ak}):$
$$\omega_A\colon A\to \mathbb R\times \ldots \times \mathbb R, 
\quad a\mapsto (\omega_{A1}(a),\ldots,\omega_{Ak}(a)).$$
Taking into account
the the direct product is a categorical direct product for real algebras, we have that,
$$M(A_1\times \ldots \times A_k) = M(A_1) \times \ldots \times M(A_k),$$
therefore multi-near-points are tuples of near-points. In particular we have
$$M(\mathbb R \times \ldots \times \mathbb R) = M \times \ldots \times M.$$
In this way we can consider tuples of points as an special case of multi-near-points. The composition
with $\omega_A$ gives the natural projection $\pi_A\colon M(A)\to M^k$.

\begin{proposition}\label{Proposition2}
Let $A = A_1\times \ldots \times A_k$ and $B = B_1\times \ldots \times B_s$ 
be multi-Weil algebras of multiplicity $k$ and $s$ respectively and let
$\varphi\colon A\to B$ be a $\mathbb R$-algebra morphism. The composition with $\varphi$
induces an smooth map $\varphi_*\colon M(A)\to M(B)$. This is a morphism of bundles in the following
sense, there is a map 
$$\sigma\colon \{1,\ldots, s \}\to \{1,\ldots, k \}$$
such that for any $p^A \in M(A)$,
$$\pi_B(p^A) = (\pi_{A\sigma(1)}(p^A), \ldots, \pi_{A\sigma(s)}(p^A)).$$ 
\end{proposition}

\proofname. It is easy to check that $\varphi$ induces a morphism 
$\sigma^*\colon \mathbb R^k \to \mathbb R^s$ such that $\sigma^*(\omega_A(a)) = \omega_B(\varphi(a))$.
Then, taking into account that ${\rm Hom}_{\mathbb R-alg}(\mathbb R^{k},\mathbb R^{s})$ is
the set of maps from $\{1,\ldots,s\}$ to $\{1,\ldots,k\}$ we finish.
\hfill$\square$

\subsection{Prolongation of functions}

Let $A$ be a multi-Weil algebra and $M(A)$ the bundle of $A$-multi-near-point 
in $M$. For each function $f\in \mathcal C^\infty_M(U)$ there is a prolongation
of $f$ to a $A$-valued function $f^A$ in $M(A)$,
$$f^A\colon U(A)\to A, \quad p^A\mapsto f^A(p(A)) = p^A(f).$$
If we consider $\{a_1,\ldots,a_s\}$ a basis of $A$, then, we can decompose
$f^A$ in real components:
$$f^A = \sum_{i=1}^s f_ia_i,$$
where each $f_i$ is an smooth function in $A(U)$. The functions 
$f_1,\ldots,f_s$ are called the real components of $f$ in $U(A)$ relative
to the basis $\{a_1,\ldots,a_s\}$.

\begin{example}
Let us consider $A = \mathbb R^k$ so that $M(\mathbb R^n) = M^k$. Le us consider 
$f\in\mathcal C^\infty(M)$. We have:
$$f^A(p_1,\ldots,p_k) = (f(p_1),f(p_2),\ldots,f(p_k)),$$
and the real components in the canonical basis of $\mathbb R^k$ are,
$$f_i(p_1,\ldots,p_k) = f(p_i).$$ 
\end{example}

\begin{example}
Let us consider $f\in\mathcal C^\infty(M)$, $A = \mathbb R[[\varepsilon]]_{1,r}$, 
and $\gamma\colon(\mathbb R,0)\to M$ so that $M(A) = T^rM$ and $j_0^r(\gamma)\in T^rM$. We have:
$$f^A(j^r_0\gamma) = f(p) + \varepsilon 
\left.\frac{d}{d\varepsilon}\right|_{\varepsilon=0}f(\gamma(\varepsilon)) + \ldots +
\frac{\varepsilon^r}{r!} 
\left.\frac{d^r}{d\varepsilon^r}\right|_{\varepsilon=0}f(\gamma(\varepsilon))
$$
And the real components in the basis $\left\{1,\varepsilon,\ldots \varepsilon^r/r!\right\}$
are:
$$f^{(i)}(j_p^r\gamma) = \left.\frac{d^i}{d\varepsilon^i}\right|_{\varepsilon=0}f(\gamma(\varepsilon)).$$
In particular, for a system of coordinates $x_1,\ldots,x_n$ in $M$, their real components
are the coordinate functions $x_i^{(j)}$ in $T^rM$.
\end{example}

\subsection{Differential invariants}

Let $A$ be a (multi-)Weil algebra and $\phi\colon M\to M$ a diffeomorphism.
Let us denote by $\phi^*$ the induced automorphism of $\mathcal C^\infty(M)$ defined
as $\phi^*(f)(p) = f(\phi(p))$. The composition with $\phi^*$ naturally gives us an 
diffeomorphism of $M(A)$. This assignation compatible with the composition, and thus
the action of $G$ in $M$ naturally prolongs to $M(A)$. 

\begin{example}
For $A = \mathbb R^k$ we have $M(A) = M^k$ and 
this prolongation is simply the diagonal action of $G$ in $M^k$.
\end{example}

\begin{example}
Let $x_1,\ldots,x_n$ be a system of coordinates in an open subset of $M$, and let 
$(g_1,\ldots,g_r)$ be a system of coordinates in a neighbourhood of the identity in $G$,
the the action of $G$ in $M$ has an expression in coordinates:
$$\tilde x_1 = X_1(x,g), \dots, \tilde x_n = X_n(x,g).$$
The corresponding open subset of tangent space $T^rM$ of order $r$ is coordinated by 
$x_1,\ldots,x_n,x_1',\ldots,x_n',\ldots,x_1^{(r)},\ldots,x_n^{r}$. The action of $G$ in $T^rM$
has an expression in coordinates that is computed iteratively:
\begin{eqnarray*}
\tilde x_j &=& X_j(x,g), \\
\tilde x_j' &=& X_j'(x,x',g) = \sum_{i=1}^n \frac{\partial X_j}{\partial x_i}(x,g)x_i', \\
 &\vdots& \\
\tilde x_j^{(r)} &=& X_j^{(r)}(x,x',\ldots,x^{(r-1)},g) = 
\sum_{i=1}^n \frac{\partial X_j^{(r-1)}}{\partial x_i}(x,g,x',\ldots,x^{(r-1)})x_i'.
\end{eqnarray*}
\end{example}

Given a Weil algebra $A$, a \emph{(local) differential invariant} of type $A$ is a 
(local) invariant of the action
of $G$ in $M(A)$. In particular, a local invariant of rank $m$ and order of the action of $G$ in $M$
is a (local) invariant of the action of $G$ in $T^{(m,r)}M$. There are formal differential operators
called the \emph{total derivatives}:
$$\frac{\mathfrak d}{\mathfrak d \varepsilon_j}\colon \mathcal C^{\infty}(T^{m,r}M) \to \mathcal C^\infty(T^{m,r+1}M)$$
whose expression in coordinates is:
$$\frac{\mathfrak d}{\mathfrak d\varepsilon_j} = \sum_{i=1}^n x_{i,e_j}\frac{\partial}{\partial x_i} +
\sum_{|\alpha|>1}\sum_{i=1}^n x_{i,\alpha+e_j}\frac{\partial}{\partial x_{i,\alpha}}$$
which, for the case $m=1$ has a simpler expression:
$$\frac{\mathfrak d}{\mathfrak d\varepsilon} = 
x_1'\frac{\partial}{\partial x_1} + \ldots + x_n'\frac{\partial}{\partial x_n} + 
x_1''\frac{\partial}{\partial x_1'} + \ldots + x_n''\frac{\partial}{\partial x_n'} + \ldots$$
If $I$ is a (local) differential invariant of rank $m$ and order $r$, then its $m$ total derivatives 
are (local) differential invariants of rank $m$ and order $r+1$. The classical finiteness
theorem of Lie (given originally in \cite[p. 760]{Lie1893}) states:

\begin{theorem}[Lie]
For each $m$ there is an $r_0$ such that any local differential invariant of rank $m$ and $r$
with $r>r_0$ of the action of $G$ in $M$ is functionally dependent of differential invariants
of rank $m$ and order $r_0$ and their total derivatives. 
\end{theorem}


\section{Derivation of joint invariants}

\subsection{Twisted differential}

The following result is a simple, but important, observation. 
It is originally stated for in \cite{Weil1953} for Weil algebras, but the proof 
remains the same for multi-Weil algebras.

\begin{proposition}[Weil]\label{PWeil}
Let $M$ be an smooth manifold, 
$A$, $B$ multi-Weil algebras. There are canonical natural diffeomorphisms,
$$M(A)(B) \simeq M(A\otimes_{\mathbb R} B) \otimes M(B)(A)$$
\end{proposition}

\proofname.
Let us consider $\{a_i\}_{i\in J}$ a basis of $A$, $\{b_j\}_{j\in J}$ basis of
$B$. Then $\{a_i\times b_j\}_{I\times J}$ is a basis of $A\otimes_{\mathbb R}B$. 
Let $x_1\ldots,x_n$ be a system of coordinates in $U\subseteq M$. Then, the real components
$x_{m,i},$ defined by:
$$x_m^A = \sum_{i\in I} x_{m,i}a_i,$$
form a system of coordinates in $U(A)$ and the real componets $x_{m,i,j}$ defined by:
$$x_{m,i}^A = \sum_{j\in J} x_{m,i,j} b_j$$
form a system of coordinates in $U(A)(B)$. On the other hand, the real components
$x_{m,(i,j)}$ defined by:
$$x_m^{A\otimes_{\mathbb R}B} = \sum_{(i,j)\in I\times J} x_{m,(i,j)} a_i\otimes b_j$$
form a system of coordinates in $U(A\otimes_{\mathbb R} B)$. The desired diffeomorphism comes
from de identification of the coordinates $x_{m,i,j}$ with $x_{m,(i,j)}$.
\hfill$\square$

\begin{corollary}\label{Corollary1}
There are canonical natural diffeormorphisms:
$$M(A)^k \simeq M(A^k) \simeq M^k(A).$$
\end{corollary}

\proofname. Take $B= \mathbb R^k$ in Proposition \ref{PWeil}.
\hfill$\square$

\medskip

Let $\bar\sigma$ be an $r$-algebra morphism $\bar{\sigma}\colon A\to A^k$.
By the universal property of direct product, $\bar{\sigma}$ is a $k$-tuple
of $\bar{\sigma} = (\sigma_1,\ldots,\sigma_k)$ of morphisms from $A$ to $A$,
$$\bar\sigma(a) = (\sigma_1(a),\ldots,\sigma_k(a)).$$

By Proposition \ref{Proposition2} and Corollary \ref{Corollary1}, 
$\bar{\sigma}$ induces an smooth map 
$$\bar\sigma_*\colon M(A)\to M(A^k)\simeq M^k(A),$$
also, for a function $I$ defined in an open subset $U\subset M^k$
the $A$-prolongation of $I$ is an $A$-valued function defined in 
$U(A)\subseteq M^k(A)$. 

\begin{definition}
Let $I$ be a function defined in some open subset $U$ of $M^k$. We call
the \emph{$\bar\sigma$-twisted differential} $D_{\bar{\sigma}}$ 
of $I$ to the $A$-valued function $D_{\bar\sigma} I = I^A\circ \bar\sigma_*$.
\end{definition}

Let us denote $\Delta^k$ the diagonal 
submanifold in $M^k$, and ${\bf di}_k\colon M\to M^k$
the diagonal map. Note that the domain of definition of 
$D_{\bar{\sigma}}I$ is ${\bf di}_k^{-1}U(A)$.
Thus, the $\sigma$-twisted differential makes sense only
for functions whose domain of definition intersects $\Delta_k$.

\begin{lemma}\label{LemmaA}
Let $U\subset M$ be an open subset and $I\subset\mathcal C^\infty U$ 
be a (local) invariant of the action of $G$ in $M$. For each multi-Weil
algebra $A$, the $A$-valued function $I^A\in\mathcal C^\infty(U(A),A)$ is
an $A$-valued (local) invariant the action of $G$ in $M(A)$.
\end{lemma}

\proofname.
The assignation $M\leadsto M(A)$ is a functor. For each smooth map $\alpha_g\colon M\to M$
we have $(I\circ \alpha_g)^A = I\circ \alpha_g^A$. Hence, if $I$ is an invariant
then $I^A$ is  an invariant.\hfill$\square$

\begin{lemma}\label{LemmaB}
Let $\sigma\colon A\to B$ be a morphism of multi-Weil algebras. Let 
$\sigma_*\colon M(A) \to M(B)$ be the induced morphism between the  
spaces of multi-near-points. Let $E$ be a vector space, and  
$I\in \mathcal C^\infty(W,E)$ defined in $W\subset M(B)$ be a (local)
invariant of the prolonged action of $G$ in $M(B)$. Then $I\circ \sigma_*$
is a (local) invariant of the action of $G$ in $M(A)$. 
\end{lemma}

\proofname.
As the above lemma, since the assignation $A\leadsto M(A)$ is a functor,  
$\sigma_*$ is a $G$-equivariant map
between $M(A)$ and $M(B)$.\hfill$\square$

\begin{theorem}
If $I$ is a (local) $k$-joint invariant of the action of $G$ in $M$, then
$D_{\bar{\sigma}}I$ is an $A$-valued (local) differential invariant of
type $A$ of the action of $G$ in $M$.
\end{theorem}

\proofname. First, we apply the Lemma \ref{LemmaA} to the case of $M^k$,
and then Lemma \ref{LemmaB} to the case of $\bar{\sigma}$. \hfill$\square$

\medskip

Thus, for any basis of $A$, the real components of $D_{\sigma}I$ are
local differential invariants of the action of $G$ in $M$. In particular, 
for the case $A = \mathbb R[[\varepsilon]]_{1,r}$ the coefficient
of the twisted derivative in $\varepsilon^r/r!$ is the $r$-th $\bar \sigma$-twisted
derivative of the invariant 
$\left(\frac{\mathfrak d^r }{\mathfrak d\varepsilon^r I}\right)_{\bar\sigma}$. 

$$D_{\bar{\sigma}} I = \sum_{i=1}^r 
\frac{\varepsilon^i}{i!}
\left(\frac{\mathfrak d^i I}{\mathfrak d\varepsilon^i}\right)_{\bar\sigma}.
$$

\begin{example}\label{Endo}
For the action of ${\rm Mov}(\mathbb R^n)$ the function
$$d^2(p_1,p_2) = \sum_{i=1}^n (x_{i,1}-x_{i,2})^2$$
is a $2$-joint invariant. 
Let us consider in $\mathbb R[[\varepsilon]]_{1,1}$ the family of endomorphisms,
$$\tau_i\colon \mathbb R[[\varepsilon]]_{1,1}\to \mathbb R[[\varepsilon]]_{1,1},
\quad \varepsilon \mapsto i\varepsilon.$$
And denote $\bar \tau = (\tau_0,\tau_1)$.
The twisted derivative is computed in the following way.
We set the $(1,1)$-near-point ,
$$p(\varepsilon) = (x_1 + \varepsilon x_1',\ldots, x_n + \varepsilon x_n'),$$
then, 
$$D_{\bar\tau}(d^2)(p(\varepsilon)) = d(p(0),p(\varepsilon))^2,$$
Where this last expression is seen as a Taylor series of order one in $\varepsilon$.  
A direct computation yields:
$$D_{\bar\tau}(d^2) =  \varepsilon \sum_{i=1}^n (x_{i})'^{2},$$
and thus,
$$\left(\frac{\mathfrak d (d^2)}{\mathfrak d\varepsilon}\right)_{\bar\tau}
= \sum_{i=1}^n (x_i')^2,$$
which is the infinitesimal quadratic expression of the euclidean metric, and it is a 
differential invariant of rank $1$ and order $1$.
\end{example}

\begin{example}
We consider the action ${\rm Mov}(2, \mathbb R)$ in $\mathbb R^2$ 
and the following invariant, proportional to the oriented area $3$-joint invariant:
$$A(p_1,p_2,p_3) =  \left| \begin{array}{cc}
x_2-x_1 & x_3-x_1 \\
y_2-y_1 & y_3-y_1
\end{array}
\right| $$ 
And denote $\bar \tau = (\tau_0,\tau_1,\tau_2)$ defined as before, but
as endomorphisms of $\mathbb R[[\varepsilon]]_{1,3}$. We set in $T^3\mathbb R^2$
the near-point:
$$p(\varepsilon) = \left(x+ \varepsilon x' + \frac{\varepsilon^2}{2}x'' +
\frac{\varepsilon^3}{6}x''',
y+ \varepsilon y' + \frac{\varepsilon^2}{2}y'' +
\frac{\varepsilon^3}{6}y'''\right)$$
A direct computation of the twisted derivative yields:  
$$D_{\bar\tau} A = A(p(0),p(\varepsilon),p(2\varepsilon))  = \varepsilon^3 \left| \begin{array}{cc}
x' & x'' \\
y' & y''
\end{array}
\right|.$$
And then, 
$$\left(\frac{\mathfrak d A}{\mathfrak d\varepsilon}\right)_{\bar\tau}
= 6(x'y''-x''y'),$$
which is a differential invariant of rank 1 and order 2.
\end{example}

\begin{example} The above example can be carried out in dimension $n$. 
We consider the action ${\rm Mov}(n, \mathbb R)$ in $\mathbb R^n$ and the 
oriented volume $(n+1)$-joint invariant:
$$V_n(p_1,\ldots,p_{n-1}) =  \det(p_{2}-p_1,p_3-p_1,\ldots,p_{n+1}-p_1),$$ 
We consider $\bar \tau = (\tau_0,\tau_1,\ldots,\tau_n)$ endomorphisms of 
$\mathbb R[[\varepsilon]]_{1,{n+1}}$ defined as in Example
\ref{Endo}. We set in $T^{n+1}\mathbb R^n$
the near-point:
$$p(\varepsilon) = \left(x_1+ \varepsilon x_1' + \ldots + \frac{\varepsilon^{n+1}}{(n+1)!}x_1^{(n+1)},
\ldots.
x_n+ \varepsilon x_n' + \ldots + \frac{\varepsilon^{n+1}}{(n+1)!}x_n^{(n+1)}
\right)$$
A direct computation of the twisted derivative yields:  
$$D_{\bar\tau} V = V(p(0),p(\varepsilon),p(2\varepsilon),\ldots,p(n\varepsilon)) = 
\varepsilon^{n+1} \Lambda_n \mathcal W(x_1',\ldots,x_n')$$
where $\Lambda_n$ is a constant that depends of $n$ with value,
$$ \Lambda_n = \frac{1}{1!2! \cdots n!}\left| \begin{array}{ccccc}
1 & 1 & 1 & \ldots & 1 \\
2 & 2^2 & 2^3 &\ldots & 2^n \\
\vdots & & & & \vdots \\
n & n^2 & n^3 &  \ldots & n^n 
\end{array}
\right|.$$
And $\mathcal W$ denotes the Wronskian determinant, 
$$ \mathcal W(x_1',\ldots, x_n') = 
\left| \begin{array}{cccc}
x_1' & x_1'' & \ldots & x_n^{(n)} \\
x_2' & x_2'' &\ldots & x_2^{(n)} \\
\vdots & & & \vdots \\
x_n' & x_n''&  \ldots & x_n^{(n)} 
\end{array}
\right|,$$
and thus, 
$$
\left(\frac{\mathfrak d^{n+1} V_n}{\mathfrak d\varepsilon^{n+1}}\right)_{\bar\tau}=
(n+1)!\Lambda_n \mathcal W(x_1',\ldots,x_n')$$
which turns out to be a differential invariant of rank $1$ and order $n$ of the
action of the movements. 

\end{example}

\subsection{Joint invariants that degenerate on the diagonal}

Let us consider $I$, $J$ two functions defined in some open subset $U$ of $M$, then
for each multi-Weil algebra $A$ it is apparent than their prolongations to $U(A)$
in $\mathcal C^\infty(U(A),A)$ satisfy: 
$$(IJ)^A = I^A J^A.$$
Here, in the left
side of the equation we have the usual product of real functions, and in the right
we have the pointwise product of elements of $A$.

Let us assume now that $A$ is a Weil algebra. If 
$J$ is a non-vanishing function, we also have that for each $p^A\in U(A)$
$J^A(p^A)\not\in\mathfrak m_A$ and thus it is an invertible element. We have
then that: 
$$ \left(\frac{I}{J}\right)^A  = \frac{J^A}{I^A}.$$
However, in some cases the above expression makes sense even if $J$ is vanishing.
Let us consider $a,b$ two elements of $A$ such that $a\in  (b)$, then, there set 
element $c\in A$ such that $bc = a$ form a class in the quotient algebra 
$A/{\bf Ann}(b)$. Thus, we may interpret the fraction as:
$$\frac{a}{b} = c \in A/{\bf Ann}(b).$$
Applying the above argument to the case of a (local) $k$-joint invariant
which is expressed as a fraction $I/J$ where the denominator $J$ vanishes
along the diagonal $\Delta_k$, we have the following result.
 
\begin{theorem}\label{singular}
Let $W\subset U$ be open subsets in $M^k$ such that:
\begin{itemize}
\item[(a)] $W$ is dense in $U$.
\item[(b)] $U$ has non empty intersection with the diagonal $\Delta_k$.
\end{itemize}
Let us denote by $U' = {\bf di}_k^{-1}(U)$. 
Let $I,J$ be smooth functions defined in $U$ such that their quotient $\frac{I}{J}$ 
is well defined in $W$ and it is a (local) $k$-joint invariant. Assume that there
is an ideal $\mathfrak p$ in $A$ such that for each $p^A\in U'(A):$
\begin{itemize}
\item[(c)] $D_{\bar\sigma}I(p^A)$ is in the ideal ($D_{\bar\sigma}J(p^A)).$
\item[(d)] ${\bf Ann}(D_{\bar\sigma}J(p^A))\subseteq \mathfrak p$.
\end{itemize}
Then, $\frac{D_{\bar{\sigma}}I}{D_{\bar\sigma}J}$ is a well defined $A/\mathfrak p$-valued
function in $U(A)$, and it is a (local) $A$-differential invariant of the action
of $G$ in $M$.
\end{theorem}

\begin{example}
Let us consider the action of ${\rm Aff}(1, \mathbb R)$ in $\mathbb{RP}_1$. 
and affine ration:
$$R(x_1,x_2,x_3) = \frac{x_3-x_1}{x_2-x_1} = \frac{I(x_1,x_3)}{J(x_1,x_2)}$$ 
We again consider $\bar \tau = (\tau_0,\tau_1,\tau_2)$ the previously defined $3$-tuple 
of endomorphisms of $\mathbb R[[\varepsilon]]_{1,2}$. Let us consider $x(\varepsilon)$
a near point in $T^2\mathbb{RP}_1$, that is expressed in the affine coordinate and its
derivatives up to second order:
$$x(\varepsilon) = x+ \varepsilon x' + \frac{\varepsilon^2}{2}x''.$$
A direct computation of the twisted derivative yields:  
$$D_{\bar\tau} I(x(\varepsilon)) = \varepsilon\left( 
x' + \varepsilon x'' + \frac{8\varepsilon^2}{6}x'''
\right),\quad
D_{\bar\tau} J(x(\varepsilon)) = \varepsilon\left( 
x' + \frac{\varepsilon}{2} x'' + \frac{\varepsilon^2}{6}x'''
\right).
$$
We are under the hypothesis of Theorem \ref{singular}, and the quotient is defined
in $\mathbb R[[\varepsilon]]_{1,2}/{\bf Ann}(\varepsilon) = \mathbb R[[\varepsilon]]_{1,1}$.
In such Weil algebra we have: 
$$\frac{D_{\bar \tau} I}{D_{\bar \tau} J} = \frac{x' + \varepsilon x''}{x' + \frac{\varepsilon}{2}x''}
= \left(x' + \varepsilon x'' \right)\left(\frac{1}{x'} - \frac{\varepsilon x''}{2 x'^2}\right)
= 1 + \frac{x''}{x'}\frac{\varepsilon}{2}.$$
It turns out that,
$$2\left(\frac{\mathfrak d^2 R}{\mathfrak d\varepsilon^2}\right)_{\bar\tau}
= \frac{x''}{x'},$$
is the logaritmic derivative, a differential invariant of rank 1 and order 2.
\end{example}

\begin{example}
Let us consider the action of ${\rm PGL}(2, \mathbb R)$ in $\mathbb{RP}_1$. 
and anharmonic ratio:
$$R(x_1,x_2,x_3) = \frac{(x_1-x_3)(x_2-x_4)}{(x_1-x_2)(x_3-x_4)} = 
\frac{I(x_1,x_3,x_3,x_4)}{J(x_1,x_2,x_3,x_4)}$$ 
We take again $\bar \tau = (\tau_0,\tau_1,\tau_2,\tau_4)$ the previously defined $4$-tuple 
of endomorphisms of $\mathbb R[[\varepsilon]]_{1,4}$. Let us consider $x(\varepsilon)$
a near point in $T^4\mathbb{RP}_1$, that is expressed in the affine coordinate and its
derivatives up to second order:
$$x(\varepsilon) = x+ \varepsilon x' + \frac{\varepsilon^2}{2}x'' + \frac{\varepsilon^3}{6}x''' +
\frac{\varepsilon^4}{24}x^{(4)}.$$
A direct computation of the twisted derivative yields:  
$$D_{\bar\tau} I(x(\varepsilon)) = I(x(0),x(\varepsilon),x(2\varepsilon),x(3\varepsilon)) = $$
$$= \varepsilon^2 \left( 4x'^2 + 12x'x''\varepsilon + \frac{24x''^2+34x'x'''}{3}\varepsilon^2 + 
o_1\varepsilon^3 \right) $$

$$D_{\bar\tau} J(x(\varepsilon)) = J(x(0),x(\varepsilon),x(2\varepsilon),x(3\varepsilon)) = $$
$$= \varepsilon^2 \left(x'^2 + 3x'x''\varepsilon + \frac{15x''^2+40x'x'''}{12}\varepsilon^2 + 
o_2\epsilon^3 \right) $$
Were $o_1$ and $o_2$ are higher order terms that do not have effect in our further computations.
We are again  under the hypothesis of Theorem \ref{singular}, and the quotient is defined
in $\mathbb R[[\varepsilon]]_{1,4}/{\bf Ann}(\varepsilon^2) = \mathbb R[[\varepsilon]]_{1,2}$.
In such Weil algebra we have: 
$$\frac{D_{\bar \tau} I}{D_{\bar \tau} J} = 
\frac{ 4x'^2 + 12x'x''\varepsilon + \frac{24x''^2+34x'x'''}{3}\varepsilon^2}
{x'^2 + 3x'x''\varepsilon + \frac{15x''^2+40x'x'''}{12}\varepsilon^2}
=$$
$$= 4 + \frac{3x''^2-2x'''x'}{x'^2} \varepsilon^2.$$
It turns out that,
$$\frac{1}{2}\left(\frac{\mathfrak d^2 R}{\mathfrak d\varepsilon^2}\right)_{\bar\tau}
= \frac{3x''^2-2x'''x'}{x'^2},$$
is the Schwartzian derivative, a differential invariant of rank 1 and order 3.
\end{example}

\subsection*{Final Remarks}
In this work we have developed part of the theory of joint invariants and explored 
the link between joint and differential invariants. There is several interesting 
open questions. One
is the theory of rational joint invariants. It is true that for an algebraic actions
the joint invariants of many points are rational functions of the liftings of a 
system of joint invariants? Another interesting question concerns to the actions of pseudogrops.
Is there also a finiteness theorem? And finally to differential invariants. It is
true that we can obtain a generating system of differential invariants by
twisted derivatives of joint invariants?

\subsection*{Acnowledgements}

Our interest about the link between differential and joint invariants 
started in a seminar with M. Malakhaltsev C. Sanabria in Universidad de los Andes. 
We want to thank them for the interesting and stimulating discussions we
carried out. We also want to thank V. Ovsienko for letting us know about
the link between the anharmonic ratio and Schwartzian derivative, and P. J.
Olver for providing interesting references. Finally  
we want to acknowledge the support of Universidad Nacional de Colombia through
research grant \emph{M\'etodos
Algebraicos en Sistemas Din\'amicos 2014} (Hermes 18636). Juan Sebasti\'an D\'iaz acknowledges
the support of Colciencias grant of the program \emph{J\'ovenes Investigadores} (Hermes 21466).

\bigskip

{\sc\noindent David Bl\'azquez-Sanz \\
Universidad Nacional de Colombia - Sede Medell\'in\\
}
E-mail: {\tt dblazquezs@unal.edu.co}

\medskip

{\sc\noindent Juan Sebasti\'an D\'iaz Arboleda \\
Universidad Nacional de Colombia - Sede Medell\'in\\
}
E-mail: {\tt jsdiaz@unal.edu.co}

\end{document}